\documentclass{amsart} \usepackage{amsmath} \usepackage{amsfonts}
\usepackage{mathrsfs} \usepackage{amsthm} \usepackage{amssymb}
\usepackage{bbm} \usepackage{bm}

\usepackage[normalem]{ulem} \usepackage{enumerate}

\usepackage[usenames]{color}

\newcommand{\bpf}[1][Proof]{{\noindent {\sc #1: }}}
\newcommand{\epf}{{{\hfill $\Box$ \smallskip}}}

 \newcommand{\good}{\mathcal{G}}
\newcommand{\R}{\mathbb{R}}
\newcommand{\DF}{F} \newcommand{\N}{\mathbb{N}}
\newcommand{\T}{\mathbb{T}} \newcommand{\Z}{\mathbb{Z}}

\newcommand{\Pp}{\mathsf{P}} 
 \newcommand{\E}{\mathbf{E}}

\newcommand{\D}{\nabla}

\newtheorem{theorem}{Theorem} \newtheorem{lemma}{Lemma}
\newtheorem{remark}{Remark} \newtheorem{corollary}{Corollary}

\title[Smooth invariant densities for random switching]{Smooth invariant densities for random switching on the torus.}

\author{Yuri Bakhtin, Tobias Hurth, Sean D. Lawley, Jonathan C.
  Mattingly} \address{Courant Institute of Mathematical Sciences, New
  York University, 251 Mercer St, New York, NY 10012 USA}
\address{Ecole Polytechnique F\'ed\'erale de Lausanne SB MATH PRST, MA B1, Station 8, CH-1015 Lausanne}
\address{Department of Mathematics, University of Utah, Salt Lake
  City, UT 84112 USA} 
\address{Mathematics Department, Duke
  University, Durham, NC 27708 USA}

\begin{document}
\begin{abstract}
  We consider a random dynamical system obtained by
 switching between the flows generated by two
  smooth vector fields on the 2d-torus, with the random switchings happening according to a Poisson process.
  Assuming that the driving vector fields are transversal to each
  other at all points of the torus and that each of them allows for a
  smooth invariant density and no periodic orbits, we prove that the
  switched system
  also has a smooth invariant density, for every switching rate. Our
  approach is based on an integration by parts formula inspired by
  techniques from Malliavin calculus. 
\end{abstract}

\maketitle

\section{Introduction}
The main goal of this paper is to prove smoothness of invariant
densities for a  class of dynamical systems generated by  random
switching between two deterministic flows defined on the two dimensional torus~$\T^2$. 
The individual flows are assumed to have everywhere positive smooth invariant densities with respect to
Lebesgue measure and have no periodic orbits or fixed points.

Many authors have studied systems with random switchings (or, {\it
  switching systems}), and they are known independently under various
titles: {\it hybrid systems}~\cite{Yin}, {\it piecewise deterministic
  Markov processes (PDMP)} (e.g.~\cite{Davis},~\cite{Malrieu_2015}),
{\it random evolutions}, see~\cite{Hersh-story:MR1962927} for the
history of the subject and extensive bibliography. Much of this work
was inpired by~\cite{Kac-telegraph:MR0510166} (a reprint of an article
published in 1956) where the first probabilistic representation of a
second-order hyperbolic equation was obtained.

The random dynamics in question can be informally described as
follows: given finitely many smooth vector fields on a manifold, a point on
the manifold follows one of them for a while and then, at a random
time, switches to one of the other vector fields chosen at random,
follows that vector field for a random time, switches again, and so
on.  We assume the times between consecutive switches are exponential with rates only depending on the current driving vector field, and independent conditioned on the sequence of driving vector fields. If the switches follow a Markov chain on the collection of vector fields, then the two-component process composed of the point on the manifold
and the driving vector field is also Markov. More general settings are
possible, for instance it is often assumed in the literature that the
rate at which switches occur depends on the point on the manifold.

Recently, several authors studied invariant measures for such
two-component Markov processes,
e.g.~\cite{Malrieu},~\cite{Faggionato},~\cite{Cloez},~\cite{LawleyMattinglyReed2015}
and~\cite{Colonius}.  Existence of an invariant measure holds true if
the manifold is compact due to a Krylov--Bogolyubov type argument
(see~\cite{Malrieu}), and can also be derived for some other systems
from similar compactness arguments or via sufficient contractivity,
see~ \cite{LawleyMattinglyReed2015}.

In~\cite{Bakhtin}, it was shown that uniqueness of an invariant
distribution follows from existence of a point $x$ that (i) is
accessible from every other point via orbits of the driving vector
fields (which can be viewed as admissible controls) and (ii) satisfies
a H\"ormander-type hypoellipticity condition, i.e., the Lie algebra
generated by the driving vector fields at $x$ coincides with the
tangent space at $x$. The same conditions guarantee the absolute
continuity of the invariant measure with respect to the volume on the
manifold. Similar results were independently obtained
in~\cite{Malrieu}, where it was also shown that, under these
assumptions, the invariant distribution is exponentially attracting in
total variation for the action of the associated Markov semigroup.

The results of~\cite{Bakhtin} and~\cite{Malrieu} can be viewed as a
simple way to look at hypoellipticity from the probabilistic
perspective. 
It is widely known that H\"ormander's hypoellipticity
conditions lead to smoothness of solutions of associated parabolic
equations, and these results have a probabilistic interpretation via
Malliavin calculus and smoothness of transition (or invariant)
densities for hypoelliptic diffusions. However, the smoothness of
invariant densities guaranteed by hypoellipticity in the diffusion
case does not hold in general for systems with switching. Even in the
simplest one-dimensional examples, invariant densities and their
derivatives may develop singularities at stable critical points of the
driving vector fields. The dynamical point of view of this phenomenon
is based on mass accumulation near criticalities due to the
exponential contraction exhibited by the flow. In~\cite{Mattingly}, an
exhaustive analysis of all kinds of singularities emerging in the
one-dimensional setting is carried out, and it is also shown that away
from the critical points the densities are smooth. The smoothness
argument is based on the fact that time averaging along an orbit of a
vector field acts as a smoothing operator.

The situation becomes more involved in higher dimensions where
singularities of the density can be created by contraction towards
attractors (potentially with complicated structure) of one vector
field and then propagated by other vector fields along their orbits.
On top of that, if one only requires hypoellipticity, the smoothing
properties of the integral operators involved are not as pronounced
and harder to exploit. In~\cite{LawleyMattinglyReed2015}, sufficient
contractivity of the system is leveraged to prove existence and
uniqueness of the invariant measure even in infinite dimensions.

In this paper, we introduce the simplest setting on the
two-dimensional torus that is devoid of the aforementioned
difficulties. Namely, we will assume that there are two driving vector
fields that are transversal to each other everywhere on the torus,
which can be interpreted as a uniform ellipticity condition. In
addition, we will impose a requirement that precludes exponential
contraction to avoid abnormal mass accumulation. Our main result is
that under these conditions, the invariant density belongs
to~$C^\infty$ (see Theorem~\ref{thm:main} at the end of
Section~\ref{sec:switching_system}).

At the core of the proof is a study of regularizing properties of the
transfer operator associated to the switching system at the moment of
time when a second switch has just occurred (see \eqref{eq:defQ1}).
Our approach is based on integration by parts with respect to times
between consecutive switches to transfer the variation in the initial point
to a variation in the noise directions. This variation in the noise
directions can then be shifted by integrating by parts to the exponential density generating
the switching times. 

This approach is inspired by the integration by
parts at the heart of Malliavin calculus which was developed initially
precisely to prove smoothness of the transition laws for stochastic differential equations
driven by white noise. There are many conceptually related works. In
\cite{BassCranston86,Bally07}, Malliavin calculus and the associated
integration by parts is developed for equations with jumps. While the
setting is different, there are some conceptual similarities. Closer
to the setting of this article \cite{Loecherbach} uses integration by
parts to study regularity of the one-dimensional marginals of the
invariant density for a class of piecewise deterministic Markov
processes with jumps. 

It is possible to study two-dimensional and higher-dimensional systems
based on vector fields that admit critical points, cycles, or
hypoellipticity points. Our progress on those systems will be reported
in another paper where we develop more delicate versions of the
methods used in the present one. Furthermore, it is relatively
straightforward to transfer the ideas here to the semi-Markov setting
when the switching times are not exponentially distributed as long as
they are given by a smooth density that decays sufficiently fast at
infinity. The ease of transferring to the semi-Markovian setting stems
largely from the fact that we work with the chain obtained after two
successive jump times essentially as in
\cite{LawleyMattinglyReed2015}.

We close the introduction with an outline of the paper. In
Section~\ref{sec:switching_system}, we describe the class of switching
systems we consider and state our main result, Theorem~\ref{thm:main}.
The proof of Theorem~\ref{thm:main} is a direct consequence of the smoothing result given in Theorem~\ref{thm:ClassicallySmooth} obtained via integration by parts, all of which is proven
in Section~\ref{sec:integration_by_parts}. In Sections~\ref{sec:aux}
and~\ref{sec:integral_equation}, we record several auxiliary
statements, notably a growth estimate on the flows generated by the
two vector fields (Lemma~\ref{lem:growth-of-jacobian}) and an integral
equation for the invariant density (Lemma~\ref{lm:integral_formula}).
This integral equation is a prerequisite for the particular
integration-by-parts argument in
Section~\ref{sec:integration_by_parts}.

{\bf Acknowledgments.} Yuri Bakhtin, Sean Lawley, and Jonathan Mattingly
gratefully acknowledge partial support from NSF via awards DMS-1460595,
DMS-RTG-1148230 and DMS-1612898 respectively.
 
\section{The switching system} \label{sec:switching_system}

We consider a switching system intermittently driven by two smooth
vector fields $u_0$ and $u_1$ on the two-dimensional torus
$\T^2 = \R^2 / \Z^2$ . Throughout the paper, smoothness means
$C^\infty$ smoothness although our results have versions involving
lower regularity requirements and conclusions.

We usually identify $\T^2$ with $[0,1)^2$ or work with the universal
cover $\R^2$. In particular, this allows us to talk about the Lebesgue
measure on the torus, and $\R^2$-vectors can serve as differences
between points on $\T^2$.


The smoothness of  $u_0$ and $u_1$ implies that for $i \in \{0,1\}$ and $x \in \T^2$, the initial-value problem 
\begin{align*}
\dot{x}(t) =& u_i(x(t)), \\
x(0) =& x\,,
\end{align*}
has a unique solution defined for all $t \in \R$. This lets us
associate flows $(x,t) \mapsto \Phi_0^t(x)$ and
$(x,t) \mapsto \Phi_1^t(x)$ to the vector fields $u_0$ and $u_1$. We
define a stochastic process $X = (X_t)_{t \geq 0}$ on $\T^2$ as
follows. Given $i \in \{0,1\}$ and $x \in \T^2$, the process $X$
follows the flow $t \mapsto \Phi_i^t(x)$, $t \geq 0$, for an
exponentially distributed random time $\tau$. Then, a switch from
$u_i$ to $u_{1-i}$ occurs and $X$ follows the flow
$t \to \Phi_{1-i}^t(y)$, $t \geq \tau$, where $y = \Phi_i^{\tau}(x)$
is the point on $\T^2$ where the switch happened. After another
exponentially distributed time, we switch back to the vector field
$u_i$, and so on.  For simplicity, we assume that the exponential
times between switches are i.i.d., so switching from $u_0$ to $u_1$
and from $u_1$ to $u_0$ happens with the same rate $\lambda>0$. While
the process $X$ by itself is not Markov, we obtain a Markov process
when adjoining a second stochastic process $A = (A_t)_{t \geq 0}$ on
the index set $\{0,1\}$ that records the driving vector field at any
given time. We denote the Markov semigroup of the two-component
process $(X,A)$ with state space $\T^2 \times \{0,1\}$ by
$(\Pp^t)_{t \geq 0}$ and the corresponding transition probability
measures by $\Pp^t_{x,i}$. A probability measure~$\mu$ on
$\T^2 \times \{0,1\}$ is called an invariant measure of
$(\Pp^t)_{t \geq 0}$ if
\begin{equation*}
\mu(E \times \{i\}) = \mu \Pp^t(E \times \{i\}) := \sum_{j \in \{0,1\}} \int_{\T^2} \Pp^t_{x,j}(E \times \{i\}) \ \mu(dx \times \{j\})
\end{equation*}
for any Borel set $E \subset \T^2$, $i \in \{0,1\}$ and $t \geq 0$.

To state our main result we need to introduce two conditions. We say
that a smooth vector field $u$ on $\T^2$ satisfies Condition A or the
{\it conjugacy condition} if the flow generated by $u$ has an invariant
measure with an everywhere positive, $C^\infty$ density with respect
to Lebesgue measure and no periodic or fixed points. We will see that
every such flow is 
smoothly conjugated to a flow with a simple structure.
We will clarify the structure of this simple flow, the conjugacy, and the role of this assumption in Section~\ref{sec:aux}. Here we only mention that presence of critical points or cycles may lead to
invariant density singularities, which happens even in one-dimensional situations studied
in~\cite{Mattingly}.

We say that a pair of two smooth vector fields $u$ and $v$ on $\T^2$
satisfies Condition~B if for every $x\in\T^2$, the vectors $u(x)$ and
$v(x)$ span the tangent space $T_x\T^2\cong \R^2$. We will also refer
to Condition~B as the {\it ellipticity} or {\it transversality}
condition. Often in this paper, we use $(u,v)$ to denote the
$2\times 2$ matrix composed of two vector columns $u$ and $v$. The
transversality condition may be rewritten as $\det(u(x),v(x))\ne 0$
for all points $x\in\T^2$.

Imposing the conjugacy conditions on individual vector fields $u_0$
and $u_1$ and the transversality condition on the pair $(u_0,u_1)$
defines a broad class of switched systems, see, e.g., Section 14.2
of~\cite{KH:MR1326374}. For example, one can start with two linear
flows on the torus, with distinct irrational slopes, and apply two
conjugations to them separately, using transformations that are
appropriately close to the identity map.

Since the torus is compact, the smoothness of driving vector fields
guarantees, by a standard Krylov--Bogolyubov argument, that there is at least one invariant measure for the Markov semigroup $(\Pp^t)_{t \geq 0}$.
Theorem~1 in~\cite{Bakhtin} ensures that $(\Pp^t)_{t \geq 0}$ admits a
unique invariant measure and that the invariant measure is absolutely
continuous with respect to the product of Lebesgue measure on $\T^2$
and counting measure on $\{0,1\}$ if there is a point $x\in\T^2$ that
(i)~satisfies a H\"ormander hypoellipticity condition and (ii)~is
accessible from any other point of the torus.

In our setting, every point $x$ on the torus satisfies these
requirements since (i)~our ellipticity condition implies the
H\"ormander condition for all points~$x \in \T^2$, and (ii) our
conjugacy condition guarantees that for any $x, y \in \T^2$, every
neighborhood of~$x$ is visited by the orbit emitted from $y$. We will
show the second part in Section~\ref{sec:aux}. We denote the unique
invariant measure by $\mu$, the marginals of $\mu$ by $\mu_0$ and
$\mu_1$, and their respective density functions with respect to
Lebesgue measure by $\rho_0$ and~$\rho_1$. We 
call~$\rho_0$ and~$\rho_1$ invariant densities.

The main result of the present paper is the following:
\begin{theorem}\label{thm:main} If smooth vector fields $u_0$ and
  $u_1$ each satisfy the conjugacy condition~A and if the pair
  $(u_0,u_1)$ satisfies the transversality condition B, then the
  invariant densities $\rho_0$ and $\rho_1$ admit $C^\infty$
  representatives for every switching rate $\lambda > 0$.
\end{theorem} We prove Theorem~\ref{thm:main} in
Section~\ref{sec:integration_by_parts}. 


\section{The basic idea} \label{sec:basic_idea}

The basic object of study will be the distribution of the process
right after two switches. In this way, our smoothing results can also
be applied to semi-Markov processes when the switching time
distribution is no longer exponential but rather some other
probability law on $(0,\infty)$ with moments of all finite orders and smooth density $\chi(t)$.  See also Remark~\ref{rm:semi_Markov} at the end of Section~\ref{sec:integration_by_parts}.  It is
natural to define the random map
$\Phi^{(S,T)}(x) = \left(\Phi_0^T \circ \Phi_1^S \right)(x)$ when $T$
and $S$ are independent, identically distributed random times with
density $\chi(t)$ each.

If $Z_0$ is distributed according to a law with density $h_0$ then the
density of the law of $Z_1=\Phi^{(S,T)}(Z_0)$, denoted by $h_1$, is
given by $h_1(x)=(Qh_0)(x)$ where $Q$ is the transfer operator defined
by
\begin{align*}
  (Qh)(x) =  \int_0^{\infty} \int_0^{\infty} \chi(s,t) \ J_{(s,t)}(x) \ h(\Psi^{(s,t)}(x)) \ ds \ dt. 
\end{align*}
Here, $\chi(s,t) = \chi(s)\chi(t)$ and
$\Psi^{(s,t)} = \big(\Phi^{(s,t)}\big)^{-1}$ and $\ J_{(s,t)}(x)$ is a
Jacobian associated with the inverse flow. All of this will be defined
precisely in Section~\ref{sec:integral_equation}. Following this imbedded chain obtained by observing the system after jumps was the perspective taken in  \cite{LawleyMattinglyReed2015}.
Our goal is to study the smoothing properties of $Q$.

Here we only want to outline the essence of the integration-by-parts
estimate at the core of our results. Let us assume that for any
direction $\xi \in \R^2$ we can find a corresponding direction
$\tau \in \R^2$ so that 
\begin{align}\label{eq:basicIdeaTransfer}
  \D_x\big(\, J_{(s,t)}(x) \ h(\Psi^{(s,t)}(x))\,\big)\xi = \D_{(s,t)}\big(\, J_{(s,t)}(x) \ h(\Psi^{(s,t)}(x))\,\big)\tau.  
\end{align} 
Then, at least formally, 
\begin{align*}
  \D_x\big( (Qh)(x) \big)\xi&=  \int_0^{\infty} \int_0^{\infty} \chi(s,t)  \D_x\big(\, J_{(s,t)}(x)
                              \ h(\Psi^{(s,t)}(x))\,\big)\xi \ ds \ dt \\
                            &=  \int_0^{\infty} \int_0^{\infty} \chi(s,t)  \D_{(s,t)}\big(\, J_{(s,t)}(x) \ 
                              h(\Psi^{(s,t)}(x))\,\big)\tau \ ds \ dt, 
\end{align*}
where we assume that the integrals converge.
Assuming $\chi(s,t)$ is smooth, then by integrating-by-parts the
derivative $\D_{(s,t)}$ can be moved onto the density $\chi(s,t)$ at
the price of generating a few boundary terms. However, none of the
terms will have any derivatives of the function $h$. Assuming that all
of the terms are well defined, we obtain an expression for $\D_x(Qh)$
which is well defined even if $h$ is not smooth. This can then be
parlayed into a proof that any invariant measure of the system must be
smooth. The precise version needed to prove our main result is
contained in Theorem~\ref{thm:IBP} and its extensions
Corollary~\ref{cor:GenIBP} and
Theorem~\ref{thm:ClassicallySmooth}. The latter two show how the above
argument can be extended to~$Q^n$ to demonstrate that every successive
application of $Q$ further smoothens the initial density.

%
%


\section{Estimates on the deterministic flows.}
\label{sec:aux}

Let us now clarify the conjugacy condition and show that it leads to
at most polynomial growth of various derivatives in time. To emphasize the generality of Condition A, we momentarily consider flows in the more general setting of a compact manifold $N$.  For every smooth vector field
$u$ on $N$,
 and for any $x \in N$, the initial-value problem 
\begin{align*}
\dot{x}(t) =& u(x(t)), \\
x(0) =& x,
\end{align*}
has a unique solution $\Phi_u^t(x)$ defined for all $t \in \R$. The
function $\Phi_u:\R\times N\to N$, $(t,x)\mapsto \Phi_u^t(x)$, is
called the flow generated by $u$. It is $C^\infty$, jointly in
$t\in\R$ and $x\in N$. We often treat the flow $\Phi_u$ as a family of
diffeomorphisms $\Phi^t_u: N\to N$, $t\in\R$.

Flows $\Phi_u$ and $\Phi_v$ generated by vector fields $u$ and $v$ on
manifolds $N$ and $M$ are smoothly conjugated if there is a $C^\infty$
diffeomorphism $\sigma:N\to M$ such that for all $t\in\R$ and
$x\in N$, \[ \Phi_u^t(x)=\sigma^{-1}\circ \Phi^t_v\circ \sigma(x). \]

One can prove (see Theorem 14.2.5 in~\cite{KH:MR1326374}) that every
smooth fixed-point-free flow on the torus conjugated to a flow
preserving a smooth positive density is also conjugated to a special
flow over a circle rotation under a smooth roof function. Let us
describe this special flow.

Let $S^1=\R^1/\Z^1$ be the unit circle, which we often view as the segment $[0,1)$ with identified endpoints.
Let $H: S^1\to (0,+\infty)$ be a smooth function and $\omega\in[0,1)$. Due to the smoothness of $H$, the set 
\[M=\{(r,h): r\in S^1,\ h\in[0,H(r)]\}\, /\sim,\] where $\sim$ is the
equivalence relation identifying points $(r,H(r))$ and $(r+\omega, 0)$
for all $r\in S^1$, is $C^\infty$ diffeomorphic to the torus so that
the flow $\tilde \Phi=\Phi_{\partial_h}$ associated with the
``vertical'' vector field $(0,1)=\partial_h$ is a smooth flow. Under
this special flow, every point $(r,h)\in M$ moves in the vertical
direction with constant speed 1, so that the $h$ component keeps
increasing until it reaches the value $H(r)$. Upon reaching
$(r, H(r))$, the point makes an instantaneous jump to $(r+\omega, 0)$
and from there continues moving upward with unit speed, etc.

Theorem 14.2.5 in~\cite{KH:MR1326374} implies that any flow associated
to a smooth vector field~$u$ on $\T^2$ satisfying Condition~A is
smoothly conjugated to a special flow with appropriately chosen
$\omega=\omega_u$ and $H=H_u$. Moreover, since Condition~A requires
that $u$ does not admit any periodic orbits, the number $\omega$ in
the above construction has to be irrational. Therefore, all orbits are
dense in $M$ for the special flow and in~$\T^2$ for $\Phi_u$. The
conjugacy of the flows $\Phi_u$ and $\tilde \Phi$ by the
diffeomorphism $\sigma$ can be rewritten as $\D_x\sigma(x) u(x)=(0,1)$
for all $x\in\T^2$. Here, $\D_x \sigma(x)$ is the Jacobian matrix of the
map~$\sigma$ at point~$x$. For fixed $t \in \R$, we denote the
Jacobian matrix of $x \mapsto \Phi^t(x)$ by $\D_x \Phi^t(x)$. For
nonnegative integers $n_1$ and $n_2$, we write
$\partial_1^{n_1} \partial_2^{n_2} \Phi^t(x)$ for the coordinatewise
partial derivative of $x \mapsto \Phi^t(x)$, where each coordinate of
$\Phi^t(x)$ is differentiated $n_1$ times with respect to the first
coordinate of $x$ and $n_2$ times with respect to the second
coordinate of $x$. Finally, for any $n \in \N$, we denote the
Euclidean norm on $\R^n$ by $\lvert \cdot \rvert$.


The following polynomial estimate on the Jacobian of the flow is
crucial for our analysis. It is based on conjugacy to a special flow
described above. This estimate implies that the Lyapunov exponents of the
flows we consider are equal to zero. If the Lyapunov exponents are
non-zero, one must be postive and one must be negative.  Excluding negative Lyapunov
exponents is natural as the associated contraction  often leads to
invariant densities with singularities.

\begin{lemma}\label{lem:growth-of-jacobian} For any smooth vector
  field $u$ on $\T^2$ satisfying Condition~A, there is a constant
$c>1$ and a family of constants $c_n > 0$, $n\ge 0$, such that for all $t>0$ and for all $x\in\T^2$, the flow $\Phi=\Phi_u$ satisfies the following estimates:
\begin{align}\label{eq:partial-derivative-bound}
  |\partial^{n_1}_1\partial^{n_2}_2\Phi^t(x)|
  \le& c_{n_1+n_2} (1+t)^{n_1+n_2}
\end{align}  
for $n_1,n_2\ge 0$ and $n_1+n_2\ge 1$, and 
\begin{align}
   c^{-1}\le \det \D_x \Phi^t(x)\le& c.  \label{eq:Jacobian-det-bounded}
\end{align}
\end{lemma}
The proof of Lemma~\ref{lem:growth-of-jacobian} is given in
Section~\ref{sec:proof_aux}. It relies heavily on smooth conjugacy of
$\Phi_u$ to a special flow described
above.

\section{An integral equation for invariant densities}
\label{sec:integral_equation}

We return to the setting from Section~\ref{sec:switching_system}. For notational convenience, we define the inverse flows 
\begin{equation*}
\Psi_i^t(x) = (\Phi_i^t)^{-1}(x) = \Phi_i^{-t}(x), \quad i \in \{0,1\}, \ t \in \R, \ x \in \T^2, 
\end{equation*}
and the composition  
\begin{equation*}
\Psi^{(s,t)}(x) = \left(\Psi_1^s \circ \Psi_0^t \right)(x), \quad (s,t) \in \R^2, \ x \in \T^2. 
\end{equation*}  
Furthermore, we define $\DF_i^t(x) = \D_x \Psi_i^t(x)$ and the Jacobian 
\begin{equation*}
J_{(s,t)}(x) = \det \left( \DF_1^s(\Psi_0^t x) \DF_0^t(x) \right). 
\end{equation*}
Finally, let 
\begin{equation}\label{eq:U}
U(x) = (u_1(x), u_0(x))
\end{equation}
be the matrix with columns $u_1(x)$ and $u_0(x)$.

Now, we extend the integral equation from Lemma 2 in~\cite{Mattingly}
to the case of the 2D-switching system introduced in
Section~\ref{sec:switching_system}. Instead of considering just the
latest switch, we consider the latest 2 switches leading to the
current state. To that end, we define the transfer operator
\begin{equation}\label{eq:defQ1}
Q h(x) = \int_{\R^2_{+}} \lambda^2 e^{-\lambda (s+t)} \ J_{(s,t)}(x) \ h(\Psi^{(s,t)}(x)) \ ds \ dt, \quad x \in \T^2   
\end{equation}  
for real-valued integrable functions $h$ on $\T^2$. Observe that if
$S$ and $T$ are independent exponentially distributed random variables
with parameter $\lambda$ then 
\begin{align}\label{eq:defQ2}
  Q h(x) = \E \left[ J_{(S,T)}(x)  h \big(\Psi^{(S,T)}(x) \big) \right]\,.
\end{align}

According to the following lemma, the invariant density $\rho_0$ is a fixed point of $Q$.

\begin{lemma}
  \label{lm:integral_formula} We have
  $\rho_0 = Q \rho_0$. 
\end{lemma}

\begin{remark} \rm For $x \in \T^2$, the
  term $Q \rho_0(x)$ can be interpreted
  as an average over possible histories
  of the previous two switches leading
  up to point $x$ and driving vector
  field $u_0$. 
\end{remark}
 
\bpf[Proof of Lemma~\ref{lm:integral_formula}] As in Lemma
2 in~\cite{Mattingly}, one can show that 
\begin{equation*}     \label{eq:int_fo_1}
\rho_i(x) = \int_{\R_{+}} \lambda e^{-\lambda t} \ \det \DF_i^t(x) \ \rho_{1-i}(\Psi_i^t(x)) \ dt, \quad i \in \{0,1\}.  
\end{equation*}
The lemma follows from plugging the instances of this identity for
$i=0$ and $i=1$ into one another and using the fact that the
pushforward of a function under the cumulative flow $\Phi_0^t \circ
\Phi_1^s$ is the composition of pushforwards under the individual
flows $\Phi_0^t$ and $\Phi_1^s$. \epf

\bigskip

\section{Smoothness through integration by parts}
\label{sec:integration_by_parts}

In this section,  we prove the main result on the smoothness of the invariant density (Theorem~\ref{thm:main}) using  integration by parts with respect to the times
between switches.
We begin by defining a collection of ``Good'' functions $\good$ for
which integration by parts can be performed.

We define $\good$ to be the set of all $C^{\infty}$ functions
$G:\T^2\times \R^2 \rightarrow \R$ such that the following conditions
hold. 
\begin{enumerate} 
\item There is a polynomial $p: \R^2 \rightarrow \R$ such that
  \begin{equation*} \lvert G(x,s,t)
    \rvert \leq p(s,t), \quad x \in \T^2, \ (s,t) \in \R^2_+.
  \end{equation*} 
\item For all $n \in \N$ and
  $\alpha = (\alpha_1, \ldots, \alpha_n)$ with $\alpha_l$ equal to
  $(s,t)$ or $x$, there is a polynomial $q:\R^2 \rightarrow \R$ such
  that 
  \begin{align*}
    \lvert \D_\alpha^n G(x,s,t)\xi \rvert \leq q(s,t) \prod_{l=1}^n \lvert \xi_l \rvert
  \end{align*}
  for all $x \in \T^2$, $(s,t) \in \R^2_+$ and $\xi = (\xi_1, \ldots, \xi_n) \in \R^{2n}$ with $\xi_i \in \R^2$ for $1 \leq i \leq n$. Here, $\nabla^n_{\alpha} G(x,s,t)$ denotes an $n$-fold differential of $G$ at the point $(x,s,t)$, which can be thought of as a multilinear form on the $n$-fold Cartesian product of $\R^2$ with itself.    
\end{enumerate}
These conditions are equivalent to saying that $G$ and all higher-order partial derivatives of $G$ are bounded by polynomials in $s$ and $t$. Observe that if $P$ is a polynomial in $n$ variables and if $G^{(1)}, \ldots, G^{(n)}$ are in $\good$, then $P(G^{(1)}, \ldots, G^{(n)})$ is in $\good$ as well. Furthermore, if $G$ and $H^{(1)}, \ldots, H^{(4)}$ are in $\good$, so is $G(H^{(1)}, \ldots, H^{(4)})$. Finally, if $G \in \good$, then the partial derivatives of $G$ of any order are in $\good$ as well.

The following lemma, which will be proven in
Section~\ref{sec:proof_aux}, shows that most objects of interest are
in $\good$.

\begin{lemma}\label{lm:smooth}
The components of $(x,s,t) \mapsto U(x)^{-1}$ and $(x,s,t) \mapsto
\Psi^t_i(x)$, $i \in \{0,1\}$, (both defined in
Section~\ref{sec:integral_equation}) are in $\good$. 
\end{lemma}

As an immediate corollary of Lemma~\ref{lm:smooth}, the components of
$(x,s,t) \mapsto \DF_i^t(x)$, $i \in \{0,1\}$, and the Jacobian
$J_{(s,t)}(x)$ (which were also defined in
Section~\ref{sec:integral_equation}) are in $\good$.

\subsection{Integration by parts}
\label{sec:integration-parts}

We begin with a remark on our notation for derivatives. If the
differential operator precedes a term in parentheses, e.g.
$\D_x (h(\Psi^{(s,t)} x))$, we apply the operator to the entire term,
so in the given example we would differentiate the function
$x \mapsto h(\Psi^{(s,t)} x)$. If we wish to differentiate only the
function $h$ and then evaluate the derivative at $\Psi^{(s,t)}(x)$, we
write $(\D_x h)(\Psi^{(s,t)} x)$.

Integration by parts applied to an integral over $[0, \infty)^2$
results in five terms: an integral over the interior of
$[0,\infty)^2$, two boundary terms corresponding to the coordinate
axes $s=0$ and $t=0$, and two boundary terms corresponding to
$s=\infty$ and $t = \infty$. In our setting, the terms corresponding
to $s=\infty$ and $t=\infty$ vanish. To deal with the remaining three
terms, we introduce the projections $\pi_0(s,t) = (s,t)$,
$\pi_1(s,t) = (s,0)$ and $\pi_2(s,t) = (0,t)$ for $(s,t) \in \R^2$.

\begin{theorem} \label{thm:IBP} Fix $G \in \good$. For any
  $\xi \in \R^2$, there exist 
  $G^{(0)}_\xi, G^{(1)}_{\xi}, G^{(2)}_{\xi} \in \good$ such that
  \begin{equation*}
    \E \left[ G(x,S,T) \, \D_x \big(h(\Psi^{(S,T)} x) \big) \xi \, \right] = \, \sum_{i=0}^2 \E \left[  G^{(i)}_{\xi}(x,\pi_i(S,T))  h \big(\Psi^{\pi_i(S,T)} x\big) \right] 
  \end{equation*}
  for all $C^1$ functions $h \colon \T^2 \rightarrow \R$. 
  In addition, there exists $K > 0$ (depending only on $G$) such that $\E\,
  \lvert G^{(i)}_\xi(x,\pi_i(S,T))\rvert \leq K \lvert \xi\rvert$ for all $i \in \{0,1,2\}$ and $x \in \T^2$, $\xi \in \R^2$.
\end{theorem}

\bpf Let $\xi \in \R^2$ and  $(s,t) \in \R_+^2$.
In the notation of Section~\ref{sec:integral_equation}, we have that
\begin{equation}\label{diff2}
  \begin{aligned}
      \D_{x} \Psi^{(s,t)}(x) \,\xi &=\DF_1^s(\Psi^t_0 x) \DF_0^t(x) \xi\,,\\\D_{(s,t)} \Psi^{(s,t)}(x) &= - F_1^s(\Psi^t_0x) U(\Psi^t_0x) \,.
  \end{aligned}
\end{equation}
The first equation is a straightforward application of the chain
rule. The second is obtained by using the forward derivative defined
by 
$\frac1\delta[\Psi^\delta_0\circ\Psi^t_0- \Psi^t_0]$ as $\delta
\rightarrow 0$ for the $t$ derivative and the backward derivative defined by $\frac1\delta[\Psi^s_1\circ\Psi^\delta_1- \Psi^s_1]$ as $\delta
\rightarrow 0$ for the $s$ derivative. 
By the uniform ellipticity condition, $U(x)$ is invertible for all $x$. 
Setting
\begin{align*}
  \tau_t(x) = -U(\Psi^t_0 x)^{-1}  \DF_0^t(x)\,,
\end{align*}
and combining the equations in  \eqref{diff2} produces
\begin{align}\label{eq:magic}
  \D_x \Psi^{(s,t)}(x) \xi = \D_{(s,t)} \Psi^{(s,t)}(x) \tau_t(x) \xi.
\end{align}
Hence, this choice of $\tau$ realizes the relationship promised in
\eqref{eq:basicIdeaTransfer} which transfers a variation in $x$ to one in $(s,t)$.

For any function $h:\T^2 \rightarrow \R$ which is $C^1$,
a direct calculation yields
\begin{align*}
  \D_x (h(\Psi^{(s,t)} x))\xi =  (\D_x h)(\Psi^{(s,t)}x)\, \D_{x} \Psi^{(s,t)}(x) \,\xi,
\end{align*}
which when combined with \eqref{eq:magic} produces 
\begin{align}\label{eq:h_magic}
    \D_x (h(\Psi^{(s,t)} x))\xi=\D_{(s,t)} (h(\Psi^{(s,t)}(x)))\,  \tau_t(x) \xi. 
\end{align}
Since $\D_x h$ is bounded and $G \in \good$, equation
\eqref{eq:h_magic} implies that  
\begin{align}   
  \E \Big[ \, G(x,S,T) \D_x &\big(h(\Psi^{(S,T)} x)\big)\,\xi \, \Big] =  \E \left[\, G(x,S,T)
                                                                            \D_{(s,t)} \big(h(\Psi^{(S,T)} x)\big)\, \tau_T(x) \xi \right] \notag \\
  =& \int_{\R^2_+} \lambda^2 e^{-\lambda (s+t)} G(x,s,t) \D_{(s,t)}(h(\Psi^{(s,t)} x)) \tau_t(x) \xi \ ds\,dt. \label{eq:ibp_prior}
\end{align}
After observing that the  components of $(x,s,t) \mapsto \tau_t(x)$ are in 
$\good$ by Lemma~\ref{lm:smooth}, we apply integration by parts
to~\eqref{eq:ibp_prior}. The divergence of the two-dimensional vector $\lambda^2 e^{-\lambda (s+t)} G(x,s,t) \tau_t(x)\xi$ with respect to $(s,t)$ equals $- \lambda^2 e^{-\lambda (s+t)} G^{(0)}_{\xi}(x,s,t)$, where
\begin{equation*}
  G^{(0)}_\xi(x,s,t) := G(x,s,t) \left(\lambda ( \mathbbm{1} \cdot \tau_t(x) \xi) - (e_2 \cdot \partial_t \tau_t(x) \xi) \right) - \D_{(s,t)} G(x,s,t) \tau_t(x) \xi. 
\end{equation*}
Here, $\cdot$ denotes the Euclidean inner product, $\mathbbm{1} = (1,1)^T$ and $e_i$ is the $i$th standard unit vector in $\R^2$ for $i \in \{1,2\}$. Since $h$ is bounded and since $G$ and the components of $\tau_t(x)$ are in $\good$, there is a polynomial $p$ such that 
\begin{equation*}
  \left \lvert G(x,s,t) h(\Psi^{(s,t)} x) (\tau_t(x) \xi \cdot e_i) \right \rvert \leq p(s,t), \quad x \in \T^2, \ (s,t) \in \R^2_+, \ i \in \{1,2\}. 
\end{equation*}
Thus, for $t \in \R_+$,  
\begin{equation*}
  \biggl \lvert \int_0^{\infty} \lambda^2 e^{-\lambda (s+t)} G(x,s,t) h(\Psi^{(s,t)} x) (\tau_t(x) \xi \cdot e_2) \ ds \biggr \rvert \leq  \int_0^{\infty} \lambda^2 e^{-\lambda (s+t)} p(s,t) \ ds, 
\end{equation*}
and the integral on the right tends to $0$ as $t \to \infty$. Similarly,
\begin{equation*}
  \lim_{s \to \infty} \int_0^{\infty} \lambda^2 e^{-\lambda (s+t)} G(x,s,t) h(\Psi^{(s,t)} x) (\tau_t(x) \xi \cdot e_1) \ dt = 0. 
\end{equation*} 
The integration-by-parts formula implies that the integral in~\eqref{eq:ibp_prior} equals 
\begin{multline*}
   \int_{\R^2_{+}} \lambda^2 e^{-\lambda (s+t)} G^{(0)}_{\xi}(x,s,t) h(\Psi^{(s,t)} x) \ ds\,dt \\
  -  \int_0^{\infty} \lambda^2 e^{-\lambda t} G(x,0,t) h(\Psi^{(0,t)} x) (\tau_t(x) \xi \cdot e_1) \ d t \\
  -  \int_0^{\infty} \lambda^2 e^{-\lambda s} G(x,s,0) h(\Psi^{(s,0)} x) (\tau_0(x) \xi \cdot e_2) \ ds.  
\end{multline*}     
The single integrals converge because $G \in \good$ and the double integral converges because all other integrals do. Defining
\begin{equation*}
  \begin{aligned}
    G^{(1)}_{\xi}(x,s,0) =& - \lambda  G(x,s,0) (\tau_0(x) \xi \cdot e_2), \\
    G^{(2)}_{\xi}(x,0,t) =& - \lambda G(x,0,t) (\tau_t(x) \xi \cdot
    e_1), 
\end{aligned} 
\end{equation*} 
we have 
\begin{equation*}
  \E \left[ \, G(x,S,T) \, \D_x \big(h(\Psi^{(S,T)} x)\big) \xi \, \right] = \, \sum_{i=0}^2  \E \left[ \,  {G}^{(i)}_\xi(x,\pi_i(S,T))\, h \big(\Psi^{\pi_i(S,T)} x \big)\, \right], 
\end{equation*}
and from Lemma~\ref{lm:smooth}, it follows that $G^{(i)}_\xi \in \good$ for $i \in \{0,1,2\}$. The second part of Theorem~\ref{thm:IBP} is a consequence of the fact that for $i \in \{0,1,2\}$, $G^{(i)}_{\xi}$ can be written as the dot product of $\xi$ and a vector-valued function whose components are in $\good$.  \epf

\bigskip

\subsection{$L^1$ smoothing estimates}

In this subsection, building on the integration by parts formula of the last section, we derive a formula for the derivative of $Qh(x)$ which does not involve the derivative of $h$.  This lets us bound the $L^1$ norm of $\nabla (Qh)$ in terms of the $L^1$ norm of $h$.  We denote the $L^1$ norm on $\T^2$ by $\| \cdot \|_{L^1(\T^2)}$ or just by $\| \cdot \|_{L^1}$.  For $n \in \N$ and a $C^n$ function $h: \T^2 \to \R$, let 
\begin{equation*}
\| \nabla^n_x h \|_{L^1} := \int_{\T^2} \sup_{\xi \in \R^{2n}: \lvert \xi_1 \rvert = \ldots = \lvert \xi_n \rvert = 1} \lvert \D_x^n h(x)\xi \rvert \ dx, 
\end{equation*}
where the supremum is taken over all $\xi = (\xi_1, \ldots, \xi_n)$ with $\xi_i \in \R^2$ and $\lvert \xi_i \rvert = 1$ for $1 \leq i \leq n$. We begin with a simple estimate which lets us bound various expectations with respect to the $L^1$ norm.

\begin{lemma}\label{lem:BasicL1Bound}
  Let $\mathcal{G}_0$ be a, possibly uncountable, subset of
  $\mathcal{G}$ such that there exists a single polynomial $p(s,t)$
  with $|G(x,s,t)|\leq p(s,t)$ for all $G \in \mathcal{G}_0$,
  $x \in \T^2$, and $s,t\geq 0 $. Then there exists a constant $K$ so
  that for any $i \in \{0,1,2\}$ and any $h \in L^1(\T^2)$,
  \begin{align*}
    \int_{\T^2}\sup_{G \in \mathcal{G}_0 }\E\big[\,|G(x,S,T)|\, |h( \Psi^{\pi_i(S,T)}(x))|\,\big] dx \leq K \|h\|_{L^1}.
  \end{align*}
\end{lemma}

\bpf Fix $i \in \{0,1,2\}$ and $h \in L^1(\T^2)$. Since $\lvert G(x,s,t) \rvert \leq p(s,t)$ for all $G \in \good_0$, we have 
\begin{equation} \label{eq:pf_BasicL1Bound} 
\int_{\T^2} \sup_{G \in \good_0} \E \big[ \lvert G(x,S,T) \rvert \lvert h(\Psi^{\pi_i(S,T)} x) \rvert \big] \ dx \leq \E \biggl[ p(S,T) \int_{\T^2} \lvert h (\Psi^{\pi_i(S,T)} x) \rvert \ dx \biggr].  
\end{equation} 
For fixed $s, t \geq 0$, we make the change of variables $y = \Psi^{\pi_i(s,t)}(x)$ and let $\Phi^{\pi_i(s,t)}$ denote the inverse of $\Psi^{\pi_i(s,t)}$. Since $\Psi^{\pi_i(s,t)}(\T^2)=\T^2$, the bound $|\det \D_x \Phi^{\pi_i(s,t)}(y)| \leq c$ from Lemma~\ref{lem:growth-of-jacobian} implies that the term on the right-hand side of~\eqref{eq:pf_BasicL1Bound} is less than or equal to 
\begin{equation*}
  c \E \big[p(S,T)\big] \int_{\T^2}  |h(y)|\, dy.
\end{equation*}
\epf

We now show the announced $L^1$ estimate on $\nabla (Q h)$ for $C^1$
functions $h$.

\begin{theorem}\label{thm:DiffL1Version}
  For any $\xi \in \R^2$, there exist
  $G^{(0)}_{\xi}, G^{(1)}_{\xi}, G^{(2)}_{\xi} \in \good$ such that
  \begin{equation} \label{eq:difFormula}
    \D_x (Q h(x)) \xi = \sum_{i=0}^2 \E \left[ G^{(i)}_{\xi}(x,\pi_i(S,T)) h \big(\Psi^{\pi_i(S,T)} x\big) \right]
  \end{equation}
for all $C^1$ functions $h$.  Furthermore, 
\begin{align} \label{eq:L_1_estimate_der}
  \|\D_x (Q h) \|_{L^1} \leq K \|h\|_{L^1}
\end{align}
for some $K>0$ independent of $h$.
\end{theorem}
\begin{remark}
  One can extend the above theorem to show that for any $h \in L^1(\T^2)$, $Q h(x)$ is in the Sobolev space $W^{1,1}$,  functions whose weak derivatives belong to $L^1(\T^2)$, with the weak derivative given by the right-hand side of \eqref{eq:difFormula}. The argument is given in the proof of Theorem~\ref{thm:ClassicallySmooth} where more is proven.
\end{remark}

\bpf[Proof of Theorem~\ref{thm:DiffL1Version} ]  Let us fix a $C^1$ function $h$.  Since $h$ is in $C^1(\T^2)$ and since $J_{(s,t)}(x)$ is in $\good$,~\eqref{eq:defQ2} implies that for any $\xi \in \R^2$, 
\begin{equation*}
\D_x (Q h(x)) \xi = \E \left[ \D_x J_{(S,T)}(x) \xi h \big( \Psi^{(S,T)} x \big) + J_{(S,T)}(x) \D_x \big( h ( \Psi^{(S,T)} x) \big) \xi \right].
\end{equation*}
Invoking again that $J_{(s,t)}(x)$ is in $\good$, we deduce from Theorem~\ref{thm:IBP} that there exist $G^{(0)}_{\xi}, G^{(1)}_{\xi}, G^{(2)}_{\xi} \in \good$, not depending on $h$, such that~\eqref{eq:difFormula} holds.  Moreover, each function $G^{(i)}_{\xi}$ can be written as the dot product of $\xi$ and a vector-valued function whose components are in $\good$ and do not depend on $\xi$.  Therefore, there exists a single polynomial $p(s,t)$ such that 
\begin{equation*}
\big \lvert G^{(i)}_{\xi}(x,s,t) \big \rvert \leq p(s,t) 
\end{equation*}
for all $i \in \{0,1,2\}$, $x \in \T^2$, $(s,t) \in \R^2_+$, and $\xi \in \R^2$ such that $\lvert \xi \rvert = 1$. By Lemma~\ref{lem:BasicL1Bound}, there exists $K > 0$ independent of $h$ such that 
\begin{equation*}
\| \D_x (Q h) \|_{L^1} \leq \sum_{i=0}^2 \int_{\T^2} \sup_{|\xi|=1} \E \left[ \big|G^{(i)}_{\xi}(x,\pi_i(S,T))\big|\, \big|h \big(\Psi^{\pi_i(S,T)} x\big)\big| \right] \ dx \leq K\|h\|_{L^1}.
\end{equation*}
\epf

\bigskip

\subsection{Smoothness}

We will now generalize the approach from the previous subsections in
order to show that the invariant density $\rho_0$ is $C^{\infty}$
smooth. In particular, we will show that for any positive integer $n$,
the derivative $\D_x^n (Q^n h)$ is bounded in $L^1$ by the $L^1$-norm
of $h$. We begin with a generalization of Theorem~\ref{thm:IBP}.
 
\begin{corollary} \label{cor:GenIBP} Let $n \geq 2$ and $G \in \good$.
  There exists $K > 0$ such that for any $\xi=(\xi_1, \dots,\xi_n)$
  with $\xi_i \in \R^2$ and for any $C^n$ function $h: \T^2 \to \R$,
  the term 
\begin{equation*} \E \left[ G(x,S,T) \D_x^n \big(
      h(\Psi^{(S,T)} x) \big) \xi \right] 
\end{equation*} 
can be written as a linear combination of integrals of the form
\begin{equation*} \E \left[ H^{(j)}_{\zeta}(x,\pi_j(S,T))
    \D_x^{n-1-k} \big( h(\Psi^{\pi_j(S,T)} x) \big)\eta \right],
\end{equation*}
  where $j \in \{0,1,2\}$, $0 \leq k \leq n-1$,
  $\eta \in \R^{2(n-1-k)}$ equal to a subset of $\xi$ with complement
  $\zeta \in \R^{2(k+1)}$, and $H^{(j)}_{\zeta} \in \good$ such that
  \begin{equation} \label{eq:H_j_estimate} \E \big\lvert
    H^{(j)}_{\zeta}(x,\pi_j(S,T)) \big\rvert \leq K \prod_{l=1}^{k+1}
    \lvert \zeta_l \rvert. 
  \end{equation} 
Neither the functions
  $H_{\zeta}^{(j)}$ from $\good$ nor the coefficients of the linear
  combination depend on $h$. \end{corollary}

\bpf Let $h$ be a $C^n$ function and let $\xi = (\xi_1, \ldots, \xi_n)$ with $\xi_i \in \R^2$ for $1 \leq i \leq n$. Set $\tilde \xi =(\xi_2, \ldots, \xi_n)$. Then, 
\begin{equation}     \label{eq:nabla_h_psi} 
\D_x^n \big( h(\Psi^{(s,t)} x) \big) \xi = \D_x^{n-1} \big( \D_x \big( h(\Psi^{(s,t)} x) \big) \xi_1 \big) \tilde \xi.
\end{equation}
By \eqref{eq:h_magic}, the right-hand
side of~\eqref{eq:nabla_h_psi}  can be written as
\begin{align*}
 \D_x^{n-1} \Big(\sum_{i=1}^2 \big[\big(\D_{(s,t)} \big( h(\Psi^{(s,t)} x) \big) e_i\big) \big] \big[\tau_t(x) \xi_1 \cdot e_i\big]\Big) \tilde
  \xi
\end{align*}
where $e_i$ the standard basis in $\R^2$. Using the
product rule, this derivative is a linear combination of terms of the form
\begin{equation*}
\big[\D_x^{n-1-k} \big( \D_{(s,t)} \big( h(\Psi^{(s,t)} x) \big) e_i \big)\eta\big] \big[\D_x^k \big( \tau_t(x) \xi_1 \cdot e_i \big)\tilde \zeta\,\big], 
\end{equation*}
where $0 \leq k \leq n-1$, $\eta \in \R^{2 (n-1-k)}$ equal to a subset
of $\tilde \xi$, and $\tilde \zeta \in \R^{2k}$ the complement of
$\eta$ in $\tilde \xi$. Fixing $k$, $i$ and $\eta$ and interchanging
the order of differentiation in the first term in the preceding product gives 
\begin{equation*}
\D_x^{n-1-k} \big( \D_{(s,t)} \big( h(\Psi^{(s,t)} x) \big) e_i \big) \eta = \D_{(s,t)} \big(\D_x^{n-1-k} \big( h(\Psi^{(s,t)} x) \big) \eta \big) e_i.  
\end{equation*}
Hence, if we set $\zeta = (\zeta_1, \ldots, \zeta_{k+1}) := (\xi_1, \tilde \zeta)$ and 
\begin{equation*}
H_{\zeta}(x,s,t) = G(x,s,t) \D_x^k \big( \tau_t(x) \xi_1 \cdot e_i \big) \tilde \zeta,
\end{equation*}
we can write 
\begin{multline}
 \E \left[ G(x,S,T) \D_x^{n-1-k} \big(\big[ \D_{(s,t)} \big( h(\Psi^{(S,T)} x) \big) e_i \big) \eta\big]\big[ \D_x^k \big( \tau_T(x) \xi_1 \cdot e_i \big) \tilde \zeta\,\big] \right]  \\
= \E \left[ H_{\zeta}(x,S,T) \D_{(s,t)} \big( \D_x^{n-1-k} \big( h(\Psi^{(S,T)} x) \big) \eta \big) e_i \right].  \label{eq:H_integral}  
\end{multline} 
The divergence of the two-dimensional vector $\lambda^2 e^{-\lambda (s+t)} H_{\zeta}(x,s,t) e_i$
with respect to $(s,t)$ is
$-\lambda^2 e^{-\lambda (s+t)} H^{(0)}_{\zeta}(x,s,t)$, where
\begin{equation*}
H^{(0)}_{\zeta}(x,s,t) = \lambda H_{\zeta}(x,s,t) - \D_{(s,t)} H_{\zeta}(x,s,t) e_i. 
\end{equation*} 
Similarly to the proof of Theorem~\ref{thm:IBP}, we obtain from integration by parts that the integral on the right side of~\eqref{eq:H_integral} equals 
\begin{equation*}
\sum_{j=0}^2 \E \left[ H^{(j)}_{\zeta}(x,S,T) \D_x^{n-1-k} \big( h(\Psi^{\pi_j(S,T)} x) \big) \eta \right], 
\end{equation*}
where 
\begin{align*}
H^{(1)}_{\zeta}(x,s,t) =& -\lambda H_{\zeta}(x,s,0) (e_i \cdot e_2), \\
H^{(2)}_{\zeta}(x,s,t) =& - \lambda H_{\zeta}(x,0,t) (e_i \cdot e_1).
\end{align*} 
Since $H_{\zeta} \in \good$, we also have $H^{(j)}_{\zeta} \in \good$ for $0 \leq j \leq 2$. In addition to $k$, $i$ and $\eta$, fix $j \in \{0,1,2\}$.  We show that there is a constant $K > 0$, independent of $\xi$ and $h$, such that~\eqref{eq:H_j_estimate} holds.  
This will complete the proof of Corollary~\ref{cor:GenIBP} because for given $n$ and $\xi$, there are only finitely many ways of choosing $k$, $i$, $j$ and $\eta$. 
As $G$ and the components of $\tau_t$ are in $\good$, $H^{(j)}_{\zeta}(x,s,t)$ can be written as 
\begin{equation}   \label{eq:sum_g}
\sum_{i_1 = 1}^2 \ldots \sum_{i_k =1}^2 \sum_{i_{k+1} = 1}^2 g_{i_1, \ldots, i_{k+1}}(x,s,t) \prod_{l=1}^{k+1} (\zeta_l)_{i_l}, 
\end{equation}
where $(\zeta_l)_{i_l}$ is the $i_l$-th component of $\zeta_l$ and $g_{i_1, \ldots, i_{k+1}}$ are functions in $\good$ that do not depend on $\xi$ or $h$. Thus, 
\begin{equation*}
\E \big\lvert H^{(j)}_{\zeta}(x,S,T) \big\rvert \leq \max_{i_1, \ldots, i_{k+1} \in \{1,2\}} \E \lvert g_{i_1, \ldots, i_{k+1}}(x,S,T) \rvert  \ 2^{k+1} \prod_{l=1}^{k+1} \lvert \zeta_l \rvert.
\end{equation*}
\epf

\bigskip

The next corollary generalizes the smoothing estimate in
Theorem~\ref{thm:DiffL1Version}.

\begin{corollary} \label{co:IBP_cons} Let $n$ be a positive integer.
  There exists $K_n > 0$ such that for any $C^n$ function
  $h: \T^2 \to \R$, we have 
  \begin{equation*} \| \D_x^n (Qh) \|_{L^1}
    \leq K_n \max\{\| h \|_{L^1}, \| \D_x h \|_{L^1}, \ldots, \|
    \D_x^{n-1} h \|_{L^1} \}. 
  \end{equation*} 
\end{corollary}

This corollary in turn implies the following result which captures the
smoothing effects of $Q^n$ and the intuition that each application of
$Q$ leads to another round of averaging; and hence, another degree of
smoothness.

\begin{corollary}\label{cor:integrationHighDer} For any $n \in \N$ there exists $K_n>0$ such that for any $C^n$ function $h: \T^2 \to \R$, we have 
  \begin{align} \label{eq:L_infty_estimate}
    \| \D_x^j (Q^n h) \|_{L^1} \leq K_n \| h \|_{L^1}, \quad 0 \leq j \leq n.
  \end{align}  
\end{corollary}

\bpf[Proof of Corollary~\ref{cor:integrationHighDer}] Applying
Corollary~\ref{co:IBP_cons} to the function $f=Q^{n-1}h$ produces
\begin{align*}
  \| \D_x^j (Q^n h) \|_{L^1} \leq C \max\{\| Q^{n-1} h \|_{L^1}, \| \D_x (Q^{n-1} h) \|_{L^1}, \ldots, \| \D_x^{j-1} (Q^{n-1}h) \|_{L^1} \}
\end{align*}
for $1 \leq j \leq n$. 
Repeatedly applying this type of estimate to each of the terms of the form $\| \D_x^{k} (Q^{n-1}h) \|_{L^1}$ on the right-hand side shows that there exists $C > 1$ so that
\begin{align*}
   \| \D_x^j (Q^n h) \|_{L^1} \leq C  \| Q^{n-j} h \|_{L^1}\,.
\end{align*}
Since $Q$ is a bounded operator on $L^1(\T^2)$, there exists $C_k$ so
that $\| Q^{k} h \|_{L^1} \leq C_k \|h \|_{L^1}$ for $0 \leq k \leq
n$. Thus, \eqref{eq:L_infty_estimate} holds with $K_n = C
\max\{C_0, \ldots, C_n\}$. \epf

\bigskip
From Corollary~\ref{cor:integrationHighDer}, we can now deduce
Theorem~\ref{thm:ClassicallySmooth} via an approximation argument. We will see that Theorem~\ref{thm:main} is essentially a corollary of this result.
\medskip
\begin{theorem}\label{thm:ClassicallySmooth}
  For any $h \in L^1(\T^2)$ and $n \in \N$, $Q^nh$ is in the Sobolev
  space  $W^{n,1}$ which consists of  functions whose weak derivatives
  up to and including order $n$ exist and are in
  $L^1(\T^2)$. Additionally $Q^{n+3}h$ is in  $C^n(\T^2)$ which is the
  space of $n$-times  continuously  differentiable functions.
\end{theorem}

\bpf
Since $h \in L^1(\T^2)$ and since $C^{\infty}(\T^2)$ is dense in $L^1(\T^2)$, there is a sequence $(h_k)_{k \geq 1}$ of $C^{\infty}$ functions that converges to $h$ in $L^1(\T^2)$. We can choose the approximating sequence in such a way that $\| h_k \|_{L^1} \leq \| h\|_{L^1}$ for all $k \geq 1$. Fix a positive integer $n$. By Corollary~\ref{cor:integrationHighDer}, we have 
\begin{equation*}
\| \D_x^j (Q^{n} h_k) \|_{L^1} \leq K_{n} \| h_k \|_{L^1} \leq K_{n} \| h \|_{L^1}
\end{equation*}
for $0 \leq j \leq n$ and $k \geq 1$. The sequence $(Q^{n} h_k)_{k \geq 1}$ is therefore a bounded sequence in the Sobolev space $W^{n,1}(\T^2)$ of $L^1$-functions whose weak derivatives  up to order $n$ are also in $L^1(\T^2)$. (As we have seen, the derivatives of $Q^{n} h_k$ even exist in the classical sense.) 

Now, the Rellich--Kondrachov theorem (see e.g.~\cite[Theorem 6.2]{Adams}) implies that $W^{n+3,1}(\T^2)$ is compactly embedded in $C^n(\T^2)$. Thus, there is a subsequence $(Q^{n+3} h_{k_i})_{i \geq 1}$ that converges to a limit in $C^n(\T^2)$. On the other hand, $(Q^{n+3} h_{k_i})_{i \geq 1}$ also converges to $Q^{n+3} h$ in $L^1(\T^2)$ because $Q$ is bounded.  This implies that $Q^{n+3} h$ has a representative in $C^n(\T^2)$.
\epf

\medskip
We now turn to the proof of the main result Theorem~\ref{thm:main} which is an immediate consequence of Theorem~\ref{thm:ClassicallySmooth} and the invariance of $\rho_0$.
\medskip

\bpf[Proof of Theorem~\ref{thm:main}] For any $n \in \N$, Theorem~\ref{thm:ClassicallySmooth} implies that  $Q^{n+3} \rho_0 \in C^n(\T^2)$. Since the invariance of $\rho_0$ implies that $\rho_0 = Q^{n+3} \rho_0$, the proof is complete. \epf

\bigskip

We now return to the proof of Corollary~\ref{co:IBP_cons}, which will
require the following lemma.

\begin{lemma} \label{lm:h_deriv_estimate} Let $n \in \Z_+$. There
  exists a polynomial $p_n(s,t)$ such that for any $C^n$ function
  $h: \T^2 \to \R$, we have 
  \begin{equation*} \| \D_x^n (h \circ
    \Psi^{(s,t)}) \|_{L^1} \leq p_n(s,t) \max\{\| h \|_{L^1}, \| \D_x
    h \|_{L^1}, \ldots, \| \D_x^n h \|_{L^1} \}. 
  \end{equation*}
\end{lemma}

\bpf We prove the lemma by induction. The base case $n=0$ follows from Lemma~\ref{lem:growth-of-jacobian} after the change of variables $y = \Psi^{(s,t)}(x)$.  In the induction step, assume that the inequality holds for some $n \in \Z_+$. For a fixed $C^{n+1}$ function $h$, $x \in \T^2$, $\xi_1 \in \R^2$ and $\xi = (\xi_2,\ldots, \xi_{n+1}) \in \R^{2n}$ with $\lvert \xi_1 \rvert = \ldots = \lvert \xi_{n+1} \rvert = 1$, we have 
\begin{equation}   \label{eq:h_derivative}
\D_x^{n+1} \big( h(\Psi^{(s,t)} x) \big) [\xi_1,\xi] = \D_x^n \big(\D_x \big(  h (\Psi^{(s,t)} x) \big) \xi_1 \big)\xi. 
\end{equation}
As 
\begin{equation*}
\D_x \big( h(\Psi^{(s,t)} x) \big) \xi_1 = (\D_x h) (\Psi^{(s,t)} x) \D_x \Psi^{(s,t)}(x) \xi_1, 
\end{equation*} by the same reasoning as in the proof of Corollary~\ref{cor:GenIBP},  
the derivative on the right side of~\eqref{eq:h_derivative} is a linear combination of terms of the form 
\begin{equation}    \label{eq:linear_comb}
 \Big[ \D_x^{n-k} \big(( \D_x h) (\Psi^{(s,t)} x) e_i \big) \eta\Big]\Big[\D_x^k \big( \D_x \Psi^{(s,t)}(x) \xi_1 \cdot e_i \big) \zeta\Big], 
\end{equation}
where $0 \leq k \leq n$, $i \in \{1,2\}$, $\eta \in \R^{2(n-k)}$ a subset of $\xi$ and $\zeta \in \R^{2k}$ the complement of $\eta$ with respect to $\xi$. For $i \in \{1,2\}$, let $g_i(y) := (\D_x h)(y) e_i$. Since $g_i$ is in $C^n$, the induction hypothesis implies that for any $k \in \{0, \ldots, n\}$ and $\eta \in \R^{2(n-k)}$, 
\begin{align}
\| \D_x^{n-k} (g_i \circ \Psi^{(s,t)}) \|_{L^1} \leq& p_{n-k}(s,t) \max\{\| g_i \|_{L^1}, \ldots, \| \D_x^{n-k} g_i \|_{L^1} \} \notag \\
\leq& p_{n-k}(s,t) \max\{\| h \|_{L^1}, \ldots, \| \D_x^{n+1} h \|_{L^1} \}.   \label{eq:p_estimate}
\end{align} 
Recall that the components of $\Psi^{(s,t)}(x)$ are in $\good$. This implies that for every $k \in \{0, \ldots, n\}$, there is a polynomial $q_k$ such that 
\begin{equation} \label{eq:q_estimate}
\lvert \D_x^k \big( \D_x \Psi^{(s,t)}(x) \xi_1 \cdot e_i \big) \zeta \rvert \leq q_k(s,t). 
\end{equation}
Here, it is important to note that the term on the right depends neither on $x$ nor on $\zeta$. Applying the estimates in~\eqref{eq:p_estimate} and~\eqref{eq:q_estimate} to the term in~\eqref{eq:linear_comb} yields the desired result. \epf

\bigskip

\bpf[Proof of Corollary~\ref{co:IBP_cons}] The case $n=1$ was treated
in Theorem~\ref{thm:DiffL1Version}, so we may assume without loss of
generality that $n \geq 2$. Let $h$ be a $C^n$ function and let
$\xi = (\xi_1, \ldots, \xi_n) \in \R^{2n}$ with
$\lvert \xi_1 \rvert = \ldots = \lvert \xi_n \rvert = 1$. Since $h$ is
assumed to be in $C^n$ and since $J_{(s,t)}(x)$ is in
$\good$,~\eqref{eq:defQ2} implies that we can write
$\D_x^n (Qh)(x)\xi$ as a linear combination of terms of the form
\begin{equation*} \E \left[ \big(\D_x^{n-k} J_{(S,T)}(x) \eta \big)\big(\D_x^k \big(
    h(\Psi^{(S,T)} x) \big) \zeta )\right], 
\end{equation*} 
where $0 \leq k \leq n$, $\eta \in \R^{2(n-k)}$ equal to a subset of $\xi$
and $\zeta \in \R^{2k}$ equal to the complement of $\eta$ in $\xi$.
Again because of $J_{(s,t)}(x) \in \good$, there are polynomials
$q_1, \ldots, q_n: \R^2 \to \R$, independent of $\xi$, such that
\begin{equation*} \left \lvert \D_x^m J_{(s,t)}(x)\eta \right \rvert
  \leq q_m(s,t) 
\end{equation*} 
for all $x \in \T^2$, $(s,t) \in \R^2_+$, $1 \leq m \leq n$, and $\eta \in \R^{2m}$ a subset
of $\xi$. By Lemma~\ref{lm:h_deriv_estimate}, there are also
polynomials $p_0, \ldots, p_n$, independent of $h$, such that
\begin{equation*} \| \D_x^m ( h \circ \Psi^{(s,t)}) \|_{L^1} \leq
  p_m(s,t) \max\{\| h \|_{L^1}, \ldots, \| \D_x^m h \|_{L^1} \}
\end{equation*} 
for $0 \leq m \leq n$. Thus, for $k < n$, we have
\begin{multline*}
  \int_{\T^2} \sup_{\xi \in \R^{2n}: \lvert \xi_1 \rvert = \ldots = \lvert \xi_n \rvert = 1} \left \lvert \E \left[ \Big(\D_x^{n-k} J_{(S,T)}(x) \eta \Big)\Big(\D_x^k \big( h(\Psi^{(S,T)} x) \big) \zeta\Big) \right] \right \rvert \ dx\\
  \leq\E \biggl[ q_{n-k}(S,T) \int_{\T^2} \sup_{\zeta \in \R^{2k}: \lvert \zeta_1 \rvert = \ldots = \lvert \zeta_k \rvert = 1} \left \lvert \D_x^k \big(h(\Psi^{(S,T)} x) \big) \zeta \right \rvert \ dx \biggr] \\
  \leq \E \big[ q_{n-k}(S,T) p_k(S,T) \big] \max\{\| h \|_{L^1},
        \ldots, \| \D_x^{n-1} h \|_{L^1}\}. 
\end{multline*} 
Moreover, we can deduce from Corollary~\ref{cor:GenIBP} that
      \begin{equation*} \E \left[J_{(S,T)}(x) \D_x^n \big(
          h(\Psi^{(S,T)} x) \big)\xi \right] 
\end{equation*} 
can be written as a linear combination of integrals of the form
\begin{equation} \label{eq:lin_comb} \E \left[
    H^{(j,k)}_{\zeta}(x, \pi_j(S,T)) \D_x^{n-1-k} \big(
    h(\Psi^{\pi_j(S,T)} x) \big)\eta \right], 
\end{equation}
where $j \in \{0,1,2\}$, $0 \leq k \leq n-1$, $\eta \in \R^{2 (n-1-k)}$ a subset of $\xi$ with complement~$\zeta$ and $H^{(j,k)}_{\zeta} \in \good$. Recall from the proof of Corollary~\ref{cor:GenIBP} that for fixed $j$, $k$ and $\eta$, $H^{(j,k)}_{\zeta}(x,\pi_j(s,t))$ can be written in the form of~\eqref{eq:sum_g}. Since the functions $g_{i_1, \ldots, i_{k+1}}$ in~\eqref{eq:sum_g} are in $\good$, there is a polynomial $q_k$, independent of $x$ and $\zeta$, such that 
\begin{equation*}
\left \lvert H^{(j,k)}_{\zeta}(x, \pi_j(s,t)) \right \rvert \leq q_k(\pi_j(s,t)).
\end{equation*}
Therefore, 
\begin{multline*}
 \int_{\T^2} \sup_{\xi \in \R^{2n}: \lvert \xi_1 \rvert = \ldots = \lvert \xi_n \rvert = 1} \left \lvert \E \left[ H^{(j,k)}_{\zeta} (x, \pi_j(S,T)) \D_x^{n-1-k} \big( h(\Psi^{\pi_j(S,T)} x ) \big)\eta \right] \right \rvert \ dx \\
\leq \max \{\| h \|_{L^1}, \ldots, \| \D_x^{n-1} h \|_{L^1} \} \E \left[ q_k(\pi_j(S,T)) p_{n-1-k}(\pi_j(S,T)) \right].
\end{multline*} 
Combining the estimates above and keeping in mind that the
coefficients in the linear combinations do not depend on $h$ or $\xi$,
we obtain the desired estimate on $\| \D_x^n (Q h) \|_{L^1}$. \epf

\bigskip

\begin{remark}    \rm    \label{rm:semi_Markov}
A close inspection of the proof of Theorem~\ref{thm:main} shows that smoothness of the invariant densities does not just hold in the case of Poissonian switching we described, but extends to semi-Markov processes for which the times between consecutive switches are distributed according to a law on $(0,\infty)$ that has a smooth density $\chi$ and admits all finite moments. Smoothness of $\chi$ is needed because we differentiate it when applying integration by parts.  The moment condition is required because, in various places, we exploit that an integral of the form 
\begin{equation*}
\int_0^{\infty} \int_0^{\infty} \chi(s) \chi(t) p(s,t) \ ds \ dt 
\end{equation*}
converges, where $p(s,t)$ is a polynomial in $s$ and $t$ that can have arbitrarily high degree. 
\end{remark}

\section{Proof of estimates on the deterministic flows} \label{sec:proof_aux}

\subsection{Proof of Lemma~\ref{lem:growth-of-jacobian}}

Let us first study the conjugated special flow $\tilde\Phi$ where the
conjugation is realized via a diffeomorphism $\sigma$. Let $t > 0$,
$x\in\T^2$ and $y=(r,h)\in M$ such that $y=\sigma(x)$. We define
$S=\{s\in(0,t]: h (\tilde \Phi^{s}(y))=0 \}$, and introduce an
ordering on $S$ by $S=\{t_1,\ldots,t_{n(x,t)}\}$ with
$t_1< \ldots < t_{n(x,t)}$. We also set $t_0=0$ and $t_{n(x,t)+1}=t$.
One can cover the trajectory $\{\tilde \Phi^s(y)\}_{s\in[0,t]}$ by a
family of $n(x,t)+1$ charts such that for $1 \leq k \leq n(x,t)+1$,
the $k$th chart contains the vertical line segment connecting
$\tilde \Phi^{t_{k-1}+0}(y)$ to $\tilde \Phi^{t_k - 0}(y)$. We can
define these charts in such a way that the flow within each chart is a
parallel translation, so in the canonical coordinates $(r,h)$ on $M$,
$\D_y \tilde \Phi^t(y)$ is the product of $n(x,t)$ Jacobian matrices
of coordinate changes between the charts. The linear map associated
with such a Jacobian matrix maps vectors $(1,H'(r_k))$ and $(0,1)$ to
$(1,0)$ and $(0,1)$, respectively,
 where $r_k=r(\tilde \Phi^{t_{k-1}}y)$.
Therefore, these matrices are given by $J_{-H'(r_k)}$, where a shear matrix $J_a$ is defined by
\[
  J_a=
  \begin{pmatrix}
    1 & 0\\
    a & 1
  \end{pmatrix},\quad a\in\R. \] Since $J_aJ_b=J_{a+b}$ for
$a,b\in\R$, we obtain that 
\begin{equation} \label{eq:D-Phi-tilde}
  \D_y \tilde \Phi^t(y)= J_{-\sum_{k=1}^{n(x,t)}H'(r_k)}=
  \begin{pmatrix}
    1 & 0\\
    -\sum_{k=1}^{n(x,t)}H'(r_k) & 1
  \end{pmatrix}. 
\end{equation} 
We immediately conclude that for all
$t$, 
\begin{equation} 
  \det \D_y \tilde\Phi^t (y)
  =1.\label{eq:Jacobian-is-1} 
\end{equation} 
Since
\begin{equation}\label{eq:Jacobian-conjugated}
  \D_x \Phi^t(x)=\nabla_x [\sigma^{-1}\circ \tilde \Phi^t \circ \sigma](x)= \nabla_y\sigma^{-1}(\tilde \Phi^t y) 
  \D_y \tilde \Phi^t(y) \nabla_x\sigma(x),
\end{equation}
we obtain due to~\eqref{eq:Jacobian-is-1}:
\begin{equation} \label{eq:Jacobian_flow}
  \det \D_x\Phi^t(x)= \det \nabla_y\sigma^{-1}(\tilde \Phi^t y) \det \nabla_x\sigma(x).
\end{equation}
The last identity
together with compactness of $\T^2$ and smoothness of $\sigma$ imply~\eqref{eq:Jacobian-det-bounded}. 
Using~\eqref{eq:D-Phi-tilde} and the identity~$(\partial_r r_k, \partial_h r_k)=(1,0)$, we obtain
that
\[
  |\partial^{n_1}_r \partial^{n_2}_h \tilde\Phi^t(y)|\le
  \Biggl|\sum_{k=1}^{n(x,t)}H^{(n_1)}(r_k)\Biggr|.
\]
Since there is $c_0>0$ such that  $n(x,t)\le c_0(1+t)$  for all $t>0$, we can use smoothness of $H$ 
and compactness of its domain to write
\begin{equation}
  \label{eq:partial-derivatives-estimates}
  |\partial^{n_1}_r \partial^{n_2}_h \tilde\Phi^t(y)|\le  c_1(1+t)  
\end{equation}
for some $c_1 > 0$ that only depends on $n_1$, and for all $t>0$.

For the remainder of the proof, we introduce the notation
$(\Phi^t_1(x), \Phi^t_2(x))$ for the coordinates of $\Phi^t(x)$ on
$\T^2$ and $(\tilde \Phi^t_r(y), \tilde \Phi^t_h(y))$ for the
coordinates of $\tilde \Phi^t(y)$ on $M$.
With~\eqref{eq:partial-derivatives-estimates} in hand, to
prove~\eqref{eq:partial-derivative-bound}, it remains to see that for
$l \in \{1,2\}$, $\partial^{n_1}_1\partial^{n_2}_2\Phi^t_l(x)$ can be
represented as a finite sum of terms of the form \[
 f(\tilde \Phi^t(\sigma(x)))g(x) \prod_{i=1}^p \partial_r^{k(i)}\partial_h^{m(i)}\tilde \Phi^t_{j(i)}(\sigma(x)),  
\]
where $f: M \to \R$ and $g: \T^2 \to \R$ are smooth functions,
$p\le n_1+n_2$, $k(i),m(i)\in\Z_+$ and $j(i)\in\{r,h\}$ for all
$i\in\{1,\ldots, p\}$. This can be checked by induction, starting
with~\eqref{eq:Jacobian-conjugated} as the induction basis.

\subsection{Proof of Lemma~\ref{lm:smooth}}

Since $U(x)^{-1}$ does not depend on $(s,t)$, we only need to verify
that $x \mapsto U(x)^{-1}$ has derivatives of all orders and that
these derivatives are bounded on $\T^2$. This follows from smoothness
of the vector fields and from the uniform ellipticity condition.

We will now show that the components of $\Psi^t_i(x)$ are in $\good$.
In this proof, we will
write the $k$th coordinate of a point $y \in \T^2$ as~$e_k \cdot y$, where
$e_1,e_2$ are the standard basis vectors inherited from $\R^2$. 

Let us fix
$i \in \{0,1\}$ and $k \in \{1,2\}$. 
As $(e_k \cdot \Psi_i^t(x))$ is bounded,
it only remains to check that its derivatives are bounded by
polynomials in $t$. For any finite sequence
$\alpha = (\alpha_1, \ldots, \alpha_n)$ of elements from $\{1,2,3\}$,
let
$\partial_{\alpha} = \partial_{\alpha_n} \partial_{\alpha_{n-1}}
\ldots \partial_{\alpha_1}$, where $\partial_1 = \partial_{x_1}$,
$\partial_2 = \partial_{x_2}$ and $\partial_3 = \partial_t$. We first
consider the case where all indices in $\alpha$ are from $\{1,2\}$,
i.e., where we only take spatial derivatives. In this situation,
Lemma~\ref{lem:growth-of-jacobian} implies that 
\begin{equation*}
  \lvert \partial_{\alpha} (e_k \cdot \Psi^t_i(x)) \rvert \leq c_n (1 + t)^n, \quad x \in \T^2, \ t \geq 0, 
\end{equation*} 
where $c_n > 0$ is some constant. The general case where
$\partial_{\alpha}$ includes a mixture of spatial and temporal
derivatives can then be reduced to the special case we have just
discussed. Namely, we will show that for any $n \in \Z_+$ and $\alpha \in
\{1,2,3\}^n$, $\partial_{\alpha} (e_k \cdot \Psi^t_i(x))$ can be
written as a polynomial in variables of the form $\partial_{\beta}
(e_l \cdot \Psi^t_i(x))$ and $(\partial_{\beta} (e_l \cdot
u_i))(\Psi^t_i x)$ for $\beta \in \bigcup_{j=0}^n \{1,2\}^j$ and $l
\in \{1,2\}$. Here, $\partial_{\alpha}$ should be interpreted as the
identity operator if $\alpha \in \{1,2,3\}^0$. This statement will
follow via a standard induction argument once we show that for $n \in
\Z_+$, $\alpha \in \{1,2\}^n$, and $m \in \{1,2,3\}$,
$\partial_m \partial_{\alpha} (e_k \cdot \Psi^t_i(x))$ and $\partial_m
((\partial_{\alpha} (e_k \cdot u_i))(\Psi^t_i x))$ can each be written
as a polynomial in variables  of the form $\partial_{\beta}(e_l \cdot \Psi^t_i(x))$ and $(\partial_{\beta}(e_l \cdot u_i))(\Psi_i^t x)$ for $\beta \in \bigcup_{j=0}^{n+1} \{1,2\}^j$ and $l \in \{1,2\}$. If $m \in \{1,2\}$, we have 
\begin{equation*} 
\partial_m \partial_{\alpha} (e_k \cdot \Psi^t_i(x)) = \partial_{\beta} (e_k \cdot \Psi^t_i(x)), \end{equation*} 
where $\beta = (\alpha, m) \in \{1,2\}^{n+1}$ is the concatenation of $\alpha$ and $m$. In addition, 
\begin{align*}
  \partial_m ((\partial_{\alpha} (e_k \cdot u_i))(\Psi^t_i x)) =& \left( \D_x \partial_{\alpha} (e_k \cdot u_i) \right)(\Psi^t_i x) \cdot \left( \partial_m \Psi^t_i(x) \right) \\
  =& \sum_{l=1}^2 (\partial_{(\alpha,l)} (e_k \cdot u_i))(\Psi^t_i x) \, \partial_m (e_l \cdot \Psi^t_i(x)), 
\end{align*}
and the right-hand side is in the desired form. 
If $m = 3$, interchanging the order of differentiation yields  
\begin{equation*}   
\partial_m \, \partial_{\alpha} (e_k \cdot \Psi^t_i(x)) = - \partial_{\alpha}(e_k \cdot u_i(\Psi^t_i x)). 
\end{equation*} 
By the chain rule for higher-order derivatives (see for instance
Theorem 2.1 in~\cite{Constantine}), the term on the right can be
written as a polynomial in variables  of the form $(\partial_{\beta} (e_k \cdot u_i))(\Psi^t_i x)$ and $\partial_{\beta} (e_l \cdot \Psi^t_i(x))$ for $\beta \in \bigcup_{j=0}^n \{1,2\}^j$ and $l \in \{1,2\}$. Finally, 
\begin{equation*}
\partial_3 ((\partial_{\alpha} (e_k \cdot u_i)) (\Psi^t_i x)) = - \sum_{l=1}^2 (\partial_{(\alpha,l)}(e_k \cdot u_i)) (\Psi^t_i x) (e_l \cdot u_i(\Psi^t_i x)). 
\end{equation*}
Since for any $\beta \in \bigcup_{j=0}^n \{1,2\}^j$ and $l \in \{1,2\}$, 
\begin{equation*}
\sup_{x \in \R^2, t \geq 0} (\partial_{\beta}(e_l \cdot u_i))(\Psi^t_i x) < \infty, 
\end{equation*}
we infer that $(e_k \cdot \Psi^t_i(x)) \in \good$.

\bibliographystyle{alpha}

\bibliography{regularity_2}

\end{document}